\documentstyle[12pt,amsfonts]{article}
\author{Kirti Joshi \& Chandrashekhar Khare}
\title{On ordinary forms \& ordinary Galois representations}
\date{}
 \begin{document}
 \maketitle
\newtheorem{thm}{Theorem}[section]              
\newtheorem{propose}{Proposition}[section] 
\newtheorem{lemma}[propose]{Lemma}
\newtheorem{cor}[propose]{Corollary}
\newtheorem{sublem}[propose]{Sublemma}
\newtheorem{schol}{Scholium}
\newtheorem{problem}{Problem}
\newtheorem{conj}{Conjecture}
\newtheorem{question}{Question}
\newtheorem{thmp}{Theorem}
\newtheorem{prop}{Proposition}
\renewcommand{\thethmp}{{\thethm${}'$}}
\newtheorem{notation}{Notation}
\newcommand{\Q}{{\bf Q}}     
\newcommand{\Z}{{\bf Z}}     
 \newcommand{\bQ}{{\bar{\Q}}}
 \newcommand{\m}{{\bf m}}
 \newcommand{\F}{{\bf F}}
 \renewcommand{\O}{{\cal O}}
 \newcommand{\bF}{{\bar{\F}}}
 \newcommand{\cond}{{\mathop{\rm cond}\nolimits}}
 \newcommand{\frob}{{\mathop{\rm Frob}\nolimits}}
 \newcommand{\tr}{{\mathop{\rm Tr}\nolimits}}
 \newcommand{\Gal}{{\mathop{\rm Gal}\nolimits}}
 \newcommand{\End}{{\mathop{\rm End}\nolimits}}

\let\isom=\simeq
\let\tensor=\otimes
\let\congruent=\equiv
\let\modulo=\bmod
\let\into=\hookrightarrow
\newenvironment{defn}{
                        \trivlist \item[\hskip \labelsep{\bf Definition}:]
                      }{
                        \endtrivlist}
\newenvironment{proof}{
                       \trivlist \item[\hskip \labelsep{\bf Proof}:]
                      }{
                        \hfill$\Box$\endtrivlist
                      }
\newenvironment{rmk}{
                        \trivlist \item[\hskip \labelsep{\bf Remark}:]
                      }{
                        \endtrivlist}

 \section{Introduction}\label{intro} 
        Let $p\geq 5$ be a prime. A well-known conjecture of Serre
(see \cite{Serre}) asserts that any two dimensional irreducible
representation
        $$ \rho: \Gal(\bQ/\Q) \to GL_2(\bF_p),$$ which is odd (in the
sense that $\det(c)=-1$, where $c$ is complex conjugation), arises
from reduction modulo $p$ of the $p$-adic representation attached to a
newform by Deligne (see \cite{Deligne}). In such a situation, we say
that $\rho$ is a {\em modular mod $p$ Galois representation}.  

To make
this more precise, let us fix once for all an embedding $\iota_p:\bQ
\into \bQ_p$. This at once provides us with a place of $\bQ$ lying
over $p$; we will denote this place by ${\frak p}$.  
The qualitative form of 
Serre's conjecture predicts that there is a newform $f$ in the
space $ S_k(\Gamma_0(N),\chi)$ of newforms of weight $k$ (throughout
this note we assume that $k\geq 2$), level $N$ and Nebentype $\chi$,
such that for any prime $\ell$ not dividing $pN$, we have
        $$ \tr(\rho(\frob_\ell)) = a_\ell \modulo {\frak p}$$ and
 $\det(\rho(\frob_\ell))=\chi(\ell)\ell^{k-1}\modulo{\frak p}$, where
 $a_\ell$ is the eigenvalue of the Hecke operator $T_\ell$, for
 $\ell\neq p$.  In fact, Serre made precise
certain minimal invariants $N(\rho), k(\rho), \varepsilon(\rho)$,
 where $\varepsilon(\rho)$ is a primitive character of conductor
 dividing $N(\rho)$, such that there is a newform $f$ of weight $k(\rho)$, level
 $N(\rho)$ and nebentype $\varepsilon(\rho)$, which gives rise to the
 representation $\rho$. We will call these the {\it Serre invariants} of
 $\rho$.  We note that $N(\rho)$ is prime to $p$, $2\leq k(\rho)\leq
 p^2-1$ and $\varepsilon(\rho)$ is a character of conductor dividing
 $N(\rho)$ and order prime to $p$; we refer the reader to \cite{Serre}
 for the precise definitions of these invariants.  As a consequence of
 the work of Mazur, Ribet, Carayol, Gross, Edixhoven and Diamond (see
 \cite{Diamond}) one now knows that if $\rho$ is modular, then it does
 indeed arise from a newform with the corresponding {\it Serre invariants}.

But there are infinitely many distinct newforms that give rise to the
same modular mod
$p$ Galois representation. We say that such forms are {\it
congruent}. We emphasise that a congruence in our sense
means that the Fourier coefficients of congruent newforms at almost
all primes are congruent mod $\frak p$. 

In this note we prove results linking the behaviour of congruent
newforms. In Section~2 we relate the mod $p$ ordinarity of a modular
mod $p$ Galois representation as above to the ordinarity (or the lack
of it) of the forms which gives rise to it (see Theorems~\ref{main-1},
\ref{main-2}, \ref{congruences}). In section 3 we state results about
the local component at $p$ of an automorphic representation and the
local representation at $p$ of the $p$-adic representation associated
to the form. These results answer what were posed as questions in
earlier versions of this article: the answers are available thanks to
\cite{Saito}.

This note was written in 1995. In spite of it not having been 
published till now, as it still gets occasionally referred to as a preprint
(\cite{Gou} for a recent instance), we felt encouraged to rescue
it from the limbo.

We thank H.~Hida, Dipendra Prasad
and J.~Tilouine for conversations.

\section{Mod $p$ representations}
        Recall that we have fixed an embedding $\iota_p:\bQ \into
\bQ_p$. This provides us with a place ${\frak p}$ lying over $p$ in $\bQ$.
Let $D_p$ be the decomposition group of $p$, with respect to the place
${\frak p}$, in $\Gal(\bQ/\Q)$; let $I_p$ be the inertia subgroup of
$D_p$. If $f$ is a newform, we will frequently use the symbol
$\pi(f)$ to denote the corresponding automorphic representation and
$\pi_p(f)$ for its $p$-component.

        We begin by recalling:

 \begin{defn} Let $\rho:\Gal(\bQ/\Q) \to GL_2(\bF_p)$ be an
 irreducible Galois representation. We say that $\rho$ is ordinary at
 $p$ if $\rho\big|_{D_p}$ is reducible and has an unramified
 one-dimensional quotient.
  \end{defn}

 \begin{defn} A
 newform of level $N$, $p$ not dividing $N$ (resp. $p$ dividing $N$)
 on $\Gamma_1(N)$ is said to be ordinary at $p$ if the eigenvalue of
 $T_p$ is a $\frak p$-adic unit (resp. the eigen value of $U_p$ on $f$ is a
 $\frak p$-adic unit). 
 \end{defn}

 \begin{rmk} If a newform $f\in S_k(\Gamma_1(N))$, with $N$ prime to $p$, is
 ordinary at $p$, then from $f$ we can create in the usual way 2
 ordinary $p$-old forms $f',f'' \in S_k(\Gamma_1(Np))$ which have the
 same eigenvalues as $f$ outside $p$ by the usual
 recipe. We recall the recipe.  If $\alpha,\beta$ are the roots of the
 Euler factor at $p$ of $L(f,s)$, then we may assume, without loss of
 generality, that $\alpha$ is a ${\frak p}$-unit. Then $f'(z)=f(z)-\beta
 f(pz)$, $f''(z)=f(z)-\alpha f(pz)$ are the required forms whose
 eigenvalues under $U_p$ are $\alpha$ and $\beta$ respectively.
\end{rmk}

 The following result of Deligne, Hida, Mazur-Wiles and Wiles (see
Section 12 of \cite{Gross} for instance), which
 is well-known, is the starting point for our note.
 
 \begin{thm}
     Let $f\in S_k(\Gamma_1(N))$ be an Hecke eigen form which is a new
form. Assume that the eigen value of $T_p$ (or $U_p$) is a $p$-adic
unit. Also assume that the mod $p$ representation $\rho_{f,p}$ is
irreducible. Then $\rho_{f,p}\big|_{D_p}$ is ordinary at
$p$.
 \end{thm}
        
This motivates the following question.
 
\begin{question}\label{main-q}
     Suppose $\rho$ is an  irreducible representation.
Assume that $\rho$ is ordinary and is modular. Then is there at least
one ordinary modular form which gives rise to $\rho$?
 \end{question}

        The following theorem provides a reasonable answer to this
question (part~1 of the following theorem can also be
found in Theorem~6.4 of \cite{Diamond}).

 \begin{thm}\label{main-1}
     Assume that we are given a modular mod $p$ Galois representation
$\rho$ as above. 
        \begin{enumerate}
            \item If $\rho$ is ordinary at $p$ then there is one (and
hence infinitely many) ordinary newform(s) which give rise to $\rho$.
        \item If $\rho$ is not ordinary, then no ordinary form can
give rise to it. 
        \item Given a $\rho$ as above, there are infinitely many
distinct newforms which are non-ordinary at $p$ which give rise to
$\rho$. Moreover these may be chosen to be either principal series or
supercuspidal at $p$. 
        \end{enumerate}
 \end{thm}

 \begin{rmk} If $f$ gives rise to $\rho$ then one can twist $f$ by a
     character of $p$-power conductor and $p$-power order to get
     non-ordinary forms which give rise to $\rho$. But using the
     theorem stated in the introduction of \cite{Khare}, we can get
     non-ordinary forms which give rise to $\rho$ in other less
     obvious ways.  
\end{rmk}

 \begin{proof}
        We assume that $\rho$ is modular and  ordinary at $p$, i.e
        $$ \rho\big|_{I_p} \sim \pmatrix{\chi & * \cr 0 & 1},$$
 where $\chi$ is  a power of the mod $p$ cyclotomic character $\omega$
(recall that $\omega:\Gal(\bQ_p/\Q_p)$ $\to$ $ \F_p^*$ is defined by
$\sigma(\zeta_p)=\zeta_p^{\omega(\sigma)}$, where $\zeta_p$ is a
primitive $p$th root of unity). Thus $\chi=\omega^a$ for some integer
$ a$. We normalise $a$ to be between $1$ and $p-1$. Then Serre defines:
        $$ k(\rho)=\cases{1+a & if $a\neq 1$\cr
                         2 & if $a=1$ and $\rho$ is peu ramifi\'ee at $p$\cr
                          p+1 & if $a=1$ and $\rho$ is tr\`es ramifi\'ee
                        at $p$}$$

        As a consequence of Edixhoven's work (see \cite{Edixhoven}),
it is known that $\rho$ arises from a form of weight $k(\rho)$ on
$\Gamma_1(N)$ for some $N$ such that $(N,p)=1$. From the above we see
that $2\leq k(\rho) \leq p+1$. We now need to recall theorems of
Deligne and Fontaine (see Theorem~2.5 and 2.6 of \cite{Edixhoven}
respectively).  
 \begin{thm}[Deligne]\label{deligne}
        Let $f$ be a newform of level $N$ prime to $p$ and weight $k$,
with $2\leq k\leq p+1$. If  $a_p\not\congruent 0\modulo{\frak p}$
then $\rho\big|_{D_p}$ is reducible and has an unramified quotient of
dimension one on which the Frobenius at $p$, ${\rm Frob}_p$, acts
by multiplication by the unique root of $x^2-a_p(f)x+p$ that is a
$\frak p$-unit and where $a_p(f)$ is the $p$th Fourier coefficient of $f$.
      \end{thm}
 \begin{thm}[Fontaine]\label{fontaine}
     With the same hypothesis as the above theorem, if $a_p\congruent
0\modulo {\frak p}$, then $\rho\big|_{D_p}$ is irreducible.
 \end{thm}

        Thus it follows from these two theorems that any  form $f\in
S_{k(\rho)}(\Gamma_1(N))$ which gives rise to $\rho$ is ordinary. 
After this the existence of infinitely many  ordinary newforms which
give rise to $\rho$ follows from the work of Hida (see \cite{Hida}).
        
        For the second assertion, we just need remark that this is a
result of Deligne, Hida and Mazur-Wiles quoted above. 

        For the third assertion, we just note that $\rho$ arises from
a newform in $S_{k(\rho)}(\Gamma_1(N(\rho))$. Then by the main theorem
of \cite{Khare} (there is a restriction on the determinant of $\rho$
in that paper which is easily seen to be superfluous in the present
application of the main theorem of that paper: \cite{Khare-Inv}
is a more recent reference), we see that $\rho$
arises from a newform in $S_2(N(\rho)p^r,\varepsilon(\rho))^{p-new}$
for all $r\geq 3$. Any such form is supercuspidal at $p$ if $r$ is odd
(see Remark~4.25 of \cite{Gelbart}). If $r$ is even then the form in
the $p$-new part of $S_2(N(\rho)p^r,\varepsilon(\rho))$ which gives
rise to $\rho$ can be chosen to be principal series at $p$.  This is
an application of Carayol's lemma (see \cite{Khare}, Remarks~11.2
and~11.3). From this the third assertion of the theorem follows from
the result that a form in the $p$-new part of
$S_2(\Gamma_0(N(\rho)p^r), \chi)$ with $\cond(\chi)=p^s$ and $s< r$
$(r>1)$ has $p$-eigenvalue zero (see \cite{Li}).  \end{proof}

  We now take up
  another aspect of Question~\ref{main-q}. We have seen in the above
  theorem that if we start with a mod $p$ modular representation
  $\rho$ which is ordinary, then any form of minimal weight and level
  prime to $p$ which gives rise to $\rho$ is ordinary.  In the theorem
  which follows we prove a similar result about the forms of weight
  two and minimal $p$-power level which give rise to $\rho$.

 \begin{thm}\label{main-2}
     Let $\rho$ be an  irreducible modular mod $p$
representation which is ordinary at $p$. Then $\rho$ arises from a form
$f\in S_2(\Gamma_1(N(\rho)p))$ with $N(\rho)$ prime to $p$.  Further
any form in $S_2(\Gamma_1(Np))$ (for any $N$ which satisfies
$(N,p)=1$) which gives rise to $\rho$ is ordinary at $p$ if
$k(\rho)\neq 2$; otherwise the form in $S_2(\Gamma_1(N))$ is ordinary
while there exists a $p$-old form in $S_2(\Gamma_1(Np))$ which is not
ordinary and which gives rise to $\rho$.
 \end{thm}
 
 \begin{proof} 

        That $\rho$ arises from $S_2(\Gamma_1(N(\rho)p))$ follows from
the level lowering results of Mazur, Ribet, Diamond and Carayol and
Edixhoven's proof of the weight part of Serre's conjecture upon using
the result (see \cite{Ash-Stevens}) that a representation $\rho$ which
arises from $S_k(\Gamma_1(N(\rho)))$ ($2\leq k \leq p+1$) also arises
from $S_2(\Gamma_1(N(\rho)p))$. 

         Suppose first that $\rho\big|_{I_p}\sim \pmatrix{\omega &
*\cr 0 & 1}$, where $\omega$ is the mod $p$ cyclotomic character. If
$\rho$ is finite at $p$ then we know that $\rho$ arises from a form
$f$ in $S_2(\Gamma_1(N))$ (for some $N$ prime to $p$, this is a
consequence of Mazur's principle cf. \cite{Ribet}) else $\rho$ arises
from a form in $S_2(\Gamma_1(N)\cap\Gamma_0(p))^{p-new}$. We see
easily from Theorem~2.5 and Theorem~2.6 of \cite{Edixhoven}, that a
newform $f$ in $S_2(\Gamma_1(N))$ which gives to $\rho$ is ordinary at
$p$; while it is easily seen that from such an $f$ we can create a
non-ordinary $p$-old form in $S_2(\Gamma_1(Np))$ which gives rise to
$\rho$.  If $\rho$ arises from
$S_2(\Gamma_1(N)\cap\Gamma_0(p))^{p-new}$ (with $(N,p)=1$) then we use
the well-known result that a form in this space has $p$-th eigen value
such that its square is $\varepsilon(p)$, where $\varepsilon$ is the
nebentype and is consequently ordinary at $p$.

        Now suppose that $\rho\big|_{I_p}\sim \pmatrix{\omega^a & *\cr
0 & 1}$ with $a\not\congruent 1\modulo{p-1}$. Then Proposition 8.13
of \cite{Gross} says that we have a Hecke equivariant
isomorphism $${\overline L}(k-2) \simeq M_k^0 \oplus
M_{p+3-k}^0[k-2]$$ for $3 \leq k \leq p$, where 
we use the notation of {\it loc. cit.} Namely, 
${\overline L}(k-2)$ is the space
of mod $p$ cusp forms on $\Gamma_1(Np)$ on which $a \in ({\bf Z}/p{\bf
Z})^*$ acts by $a^{k-2}$, i.e., 
the $k-2$nd power of the identity character, and $M_k^0$ (resp., $M_{p+3-k}^0[k-2]$) 
is the space of mod $p$ cusp forms on $\Gamma_1(N)$ of weight $k$ (resp.,
the mod $p$ space of cusp forms on $\Gamma_1(N)$ and weight $p+3-k$ 
with the Hecke action twisted by
the $k-2$nd power of the identity character as in the discussion before {\it
loc. cit.}). Now on the one hand
irreducible mod $p$ representations which arise from $M_{p+3-k}^0[k-2]$ are never
ordinary, as such representations which are reducible on restriction
to $D_p$ are of the form $$\left(\matrix{\chi&*\cr
                                    0&\chi^{k-2}}\right).$$
On the other hand as we have seen in the proof of Theorem \ref{main-1}, it follows from
Theorems \ref{deligne} and \ref{fontaine}, that a newform in 
$f \in S_k(\Gamma_1(N))$ ($(N,p)=1, 2 \leq k \leq p+1$)
gives rise to an irreducible mod $p$ representation that is ordinary
at $p$ if and only if $f$ is ordinary.

\end{proof}

 \begin{rmk}
     This result is essentially best possible as there are forms in
$S_2(\Gamma_0(Np^2))^{p-new}$ (which can be even supercuspidal at $p$:
see Remark~11.1 of \cite{Khare}) which give rise to ordinary mod
$p$ representations though any form in this space has $p^{th}$
eigenvalue zero. This may be seen from Theorem~3 of \cite{Khare},
where it is shown that if a mod $p$ representation arises from
$S_2(\Gamma_0(N))$, then it also arises from the $p$-new part of
$S_2(\Gamma_0(Np^2))$: see also \cite{Khare-Inv}.
 \end{rmk}

We prove a result using the methods here 
that responds to a question that was posed to one of
us by Dipendra Prasad 

\begin{thm}\label{congruences}
  If two newforms $f,f'$ in $S_k(\Gamma_0(N))$ with $N$ prime to $p$
  and $1 \leq k \leq p+1$ give rise to the same irreducible mod $\frak
  p$ representation, then their
  $p$th Fourier coefficients are also congruent.
\end{thm}

\begin{proof}
The weight 1 case is easy as then the mod $\frak p$ 
representations corresponding to $f,f'$ are unramified at $p$
and the $p$th Fourier coefficients are recovered as the
traces of ${\rm Frob}_p$.
So we assume that $2 \leq k \leq p+1$.
We claim that the $p$th Fourier coefficients $a_p(f)$
and $a_p(f')$ of $f$ and $f'$ are either both divisible by $\frak p$
or they are both units at $\frak p$. This claim follows from
Theorems \ref{deligne} and \ref{fontaine}. Then as $k \geq 2$
at most one of the roots of $x^2-a_p(f)x+p^{k-1}$ (resp., 
$x^2-a_p(f')x+p^{k-1}$) is a unit at $\frak p$ and if one of them
is a unit it arises as the eigenvalue of ${\rm Frob}_p$ acting
on the unique unramified quotient of the mod $\frak p$ representation
associated to $f$ (resp, $f'$). As we are assuming that the mod $\frak
p$ representations associated to $f,f'$ are isomorphic, and hence in
particular their restrictions to $D_p$ are isomorphic, we have
completed the proof.
\end{proof}

\begin{rmk}
  The above theorem is false when we do not have hypotheses
  on levels as follows from the first remark of this section.
  This issue is related to the results
  of \cite{Gross} about {\it companion forms}, and
  it is also this failure that results in the {\it twin form} phenomena of 
  \cite{Gouvea-Mazur}. We do not know what happens for weights $>p+1$
  and levels prime to $p$.
\end{rmk}

\section{$p$-adic representations}

\subsection{$p$-adic ordinarity}

Now we turn to the $p$-adic analogue of the above question. 
At the time we wrote this note in 1995 we could answer these $p$-adic
analogues in some cases using the results of Wiles. But now it is
indeed possible to give clean answers as a consequence
of the work of \cite{Saito} (see also \cite{Mokrane}). 

To fix notations, we work with a ring $\O_p$ that is the ring of
integers of a finite extension of ${\bf Q}_p$, which we will
assume to be large enough so that the eigen values of Hecke operators
of the form in question lie in it. Typically, we can take this to be the
ring of integers in the completion of the coefficient field $E_f$ of
$f$ along the valuation given by ${\frak p}$. We now recall the definition
of ordinarity.

\begin{defn}
    Let $\sigma: \Gal(\bQ/\Q) \to GL_2(\O_p)$ be an absolutely
irreducible Galois representation. We say that it is ordinary at $p$
if $\sigma\big|_{D_p}$ is reducible and has an unramified
one-dimensional quotient.
\end{defn}
 
 \begin{thm}\label{ordinary}
Let $\sigma:\Gal(\bQ/\Q)\to GL_2(\O_p)$ be the $\frak p$-adic representation
attached to a newform $f$ such that the residual representation is
absolutely irreducible and such that $\sigma$ is ordinary at $p$. Then
$f$ is ordinary.
\end{thm} 

  \begin{thm}\label{supercuspidal}
 Let $f$ be a newform such that the
corresponding automorphic representation $\pi_p(f)$ is supercuspidal
and let $\sigma$ be the $p$-adic representation associated to
$f$. Then $\sigma\big|_{D_p}$ is irreducible.
  \end{thm}

\begin{proof}
The proofs of both the theorems follow from results in \cite{Saito}. For Theorem \ref{ordinary} 
see Proposition 3.6 of \cite{Mokrane}. Theorem \ref{supercuspidal}
follows from the main theorem of \cite{Saito} and the fact that
the restriction of the Weil-Deligne parameter of a supercuspidal representation
of $GL_2({\bf Q}_p)$ to the Weil group of ${\bf Q}_p$ is irreducible.
\end{proof}

\begin{rmk} As noted above it is also 
possible to give a different proof in many cases using the
theorems of the type proven by Wiles (see \cite{Wiles-F}) identifying
ordinary deformation rings and ordinary Hecke rings of fixed tame
level, under some technical conditions. 
\end{rmk}

\subsection{Images of $p$-adic representations}

	One can also ask for more detailed information about the local
image at $p$, $\sigma|_{D_p}$, of the $p$-adic representation $\sigma$
associated to a newform $f$, in terms of $\pi_p(f)$ in general. If $f$
is ordinary, then by Deligne's theorem we have a good answer, though
we do not know if this representation can ever be semisimple upon
restriction to an open subgroup of the inertia group. We have the
following result in the weight 2 case whose argument was pointed out
to us by Dipendra Prasad:

\begin{prop} 
	Let $f$ be a non-CM newform of level prime to $p$, weight 2,
and {\it is} ordinary at $\frak p$. Then the corresponding $\frak
p$-adic representation $\sigma$ associated to $f$ is non-semisimple
when restricted to any open subgroup of $D_p$.
\end{prop}

\begin{proof} 
	This follows from the theory of Serre-Tate liftings.  Namely
from Theorem 2 of A.2.3 of Chapter IV of \cite{Serre1}, assuming the
contrary leads to the conclusion that $f$ has CM.
\end{proof}

	 The following result is a simple consequence
of Fontaine's theory (see \cite{fontaine}) but seems worth spelling
out as we find the result striking (see \cite{Volkov} for a
much more extensive study of $p$-adic representations that
arise from elliptic curves over $\Q_p$).

\begin{thm} If $E_{/{{\bf Q}_p}}$ is an elliptic curve with good 
supersingular reduction.  Let $V_p$ the representation of the Galois
group $G={\rm Gal}(\bar{\Q}_p/\Q_p)$ on $V_p=H^1_{et}(E\tensor
{\bar\Q_p},\Q_p)$. The image of $G$ in $GL(V_p)$ is dihedral, i.e.,
all supersingular elliptic curves over ${\bf Q}_p$ have formal CM.
\end{thm}

\begin{proof}
	Let $E$ be an elliptic curve with good supersingular reduction
over $\Q_p$. We choose a model for it over $\Z_p$ will write $E_0$ for
the special fibre of this model. What is being asserted here is
essentially independent of the choice of the integral model.

	Fontaine's theory provides a complete description of $V_p$ in
terms of the crystalline cohomology $H^1_{cris}(E_0/\Z_p)$ of the
special fibre $E_0$ together with a Frobenius map and a Hodge
filtration arising from comparison of $H^1_{cris}(E_0/\Z_p)$ and the
de Rham cohomology of $E/\Q_p$. Let $V_{dr}=H^1_{dR}(E/\Q_p)$,
together with its Hodge filtration $F^1H^1_{dR}(E/\Q_p)\subset
H^1_{dR}(E/\Q_p)$. Let$V_{cris}=H^1_{cris}(E_0/\Z_p)\otimes\Q_p$. This
comes equipped with a Frobenius automorphism (this is because we are
over $\Q_p$) $\phi:V_{cris}\to V_{cris}$. By a well-known theorem of
Berthelot, there is a canonical $V_{dr}\isom V_{cris}$.

	By \cite{fontaine} we know that the $G$-representation $V_p$ is completely
determined by $ V_{cris}, \phi$ and the Hodge filtration, and in
particular $\End_{G}(V_p)$ (as a Galois module) is the same as the
endomorphism algebra of the filtered $\phi$-module $V_{cris}$. So we now
compute the endomorphisms of the filtered $\phi$-module
$V_{cris}$. Because our curve has good supersingular reduction and
because we are over $\Q_p$ the characteristic polynomial of $\phi$ is
$X^2+p$, more over $\phi$ is semi-simple. Further we may assume
without loss of generality that $\phi$ is given by the matrix
	$$\phi=\pmatrix{ 0 & 1 \cr -p & 0}$$ for some basis $e_1,e_2$
of $V_{cris}$. The matrix of $\phi$ with respect to any other basis
$e_1',e_2'$ is then of the form $A\phi A^{-1}$ for a suitable $A\in
GL_2(\Q_p)$.

Note that since our curve $E_0$ is defined over $\F_p$ and is
supersingular, not all its endomorphisms are defined over $\F_p$. All
the endomorphisms of $E_0$ become visible over $\F_{p^2}$. The key
point is the following elementary fact that the endomorphisms of the
pair $(V'_{cris},\phi')$ is a quaternion algebra ramified at
$p,\infty$ given explicitly as $\pmatrix{a & b\cr pb^{\sigma} &
a^\sigma}$. From this by a simple calculation that we omit 
it follows that the endomorphisms
of the filtered module $(V_{cris},\phi)$ which are defined over 
$\Q_p$ are the scalars (we
do a more involved calculation of the same type below).

We pass to $\F_{p^2}$ to do the main calculation of endomorphisms
of filtered modules which shows that
the endomorphisms of $V_p$ ($\otimes \Q_p$) considered as a  
module for the Galois group of $\Q_{p^2}$, the unramified quadratic extension of
$\Q_p$, are $\Q_{p^2}$.  It is easy to
verify (using for instance the proper base change theorem) that the
crystalline cohomology of the curve over $\F_{p^2}$ can be obtained
from the curve over $\F_p$ as follows: the underlying
$\Q_{p^2}$-vector spaces is
$$V'_{cris}=\Q_{p^2}\tensor_{\Q_p}V_{cris},$$ we have to define the
Frobenius on $V'_{cris}$. This is done in the obvious way:
$$\phi':V'_{cris}\rightarrow V'_{cris}$$ is defined by
$$\phi'(x\tensor v)=\sigma(x)\tensor \phi(v),$$ where
$\phi:V_{cris}\rightarrow V_{cris}$ is the Frobenius on $V_{cris}$ and
$\sigma:\Q_{p^2}\rightarrow \Q_{p^2}$ is the non-trivial Galois
automorphism of $\Q_{p^2}/\Q_p$. In particular: $\phi'$ is
$\sigma$-semilinear automorphism (hence it is $\Q_p$-linear). Also
note that Fontaine's comparison theorem is compatible with this base
change as $B_{cris}$ is a $\Q_p^{nr}$-algebra.

Now one reduces to doing some elementary semilinear algebra.  We want
to compute the endomorphisms of $V'_{cris}$ as a filtered
$\phi'$-module.  An endomorphism of such data is a
linear map of $\Q_{p^2}$-vector spaces, which commutes with the
Frobenius $\phi'$ and preserve the filtration.   To compute explicitly
it is convenient do our calculations for a particular form of the
Frobenius: $$\phi:V \to V$$ given by the matrix $$\phi=\pmatrix{ 0 &
1\cr -p & 0}$$

Notice that $\phi'$ is also given by the matrix above for $\phi$ 
(merely extend the basis of
$V$ to $V'$). But now $\phi'$ is $\sigma$-semilinear: so on the
coefficients of the basis vectors $\sigma$ will operate. Let
$\tau:V'\rightarrow V'$ be any endomorphism of $V'$. We may assume
without loss of generality that $\tau$ preserves our basis vector
$e_1$. So it is upper triangular. Thus we can suppose
$\tau=\pmatrix{\alpha &
\beta\cr 0 &\delta}$, with respect to our basis. The condition that
$\phi'\tau=\tau\phi'$ means that the matrices corresponding to them
$\sigma$-commute, in symbols:
	$$ \pmatrix{ \alpha & \beta\cr 0 &\delta}\pmatrix{ 0 & 1 \cr -p
& 0}=\pmatrix{0 & 1 \cr -p & 0}\pmatrix{\alpha^\sigma &
\beta^\sigma\cr 0 &\delta^\sigma}.$$

	This gives us the relations: $\delta=\alpha^\sigma$ and
$\beta=0$. So any $\tau$ which commutes with $\phi'$ and which
preserves the basis $e_1,e_2$ must be of the form
	$$\tau=\pmatrix{ \alpha & 0 \cr 0 & \alpha^\sigma}$$ for some
$\alpha\in\Q_{p^2}$.

 	Now to recover results for our filtration basis $e_1',e_2'$ we
note that if $\tau\phi'=\phi'\tau$ for the above basis, we have
the identity
	$$ (A\tau A^{-1}) (A\phi'A^{-1})=(A\phi'A^{-1})(A\tau^\sigma
A^{-1}), $$ which says, as $A^\sigma=A$, that $A\tau A^{-1}$
$\sigma$-commutes with $A\phi'A^{-1}$ if $\tau$ $\sigma$-commutes with
$\phi'$ (with respect to the new basis).
	
So we see that if we restrict our Galois representation to $\Q_{p^2}$,
then the algebra of endomorphisms of the Tate module of our curve $E$
is isomorphic to the two dimensional $\Q_p$ vector-space
$\Q_{p^2}$. By a standard argument using the double commutant theorem
or more directly one deduces that the Lie algebra of the image of $G$
(and more generally of any open subgroup of finite index) is
toral. Hence by the irreducibility of the the representation $V_p$
that we saw earlier, we conclude that the image of $G$ is
dihedral.
\end{proof}

\begin{rmk}
	A more direct argument using formal
groups to prove this result was given by Dipendra Prasad, and is as follows.
The formal group of an elliptic curve
over ${\bf Q}_p$ with good supersingular reduction is a height 2 formal
group which is a Lubin-Tate formal group over the quadratic unramified
extension of ${\bf Q}_p$. Hence from the Lubin-Tate theory, the
maximal compact subring of this quadratic unramified extension operates
on the formal group, and thus the elliptic curve has formal CM.
\end{rmk}


\noindent{\it Addresses of the authors:} School of Mathematics, TIFR, Homi Bhabha Road,
Mumbai 400 005, INDIA. e-mail address: kirti@math.tifr.res.in, shekhar@math.tifr.res.in

 \end{document}